\documentclass[12pt]{article}

\usepackage{latexsym}
\usepackage{amsfonts}
\usepackage{graphicx}

\newtheorem{thm}{Theorem}
\newtheorem{lemma}{Lemma}
\newtheorem{ques}{Question}
\newtheorem{conj}[ques]{Conjecture}

\newcommand{\qed}{\hfill\mbox{$\framebox(5,5)[]{}$}}
\newenvironment{proof}{\par \trivlist
 \itemindent\parindent \item[\hskip\labelsep\sc Proof.]
 \ignorespaces}{\qed\endtrivlist}

\newcommand{\dist}{\mbox{dist}}
\newcommand{\ev}{\mbox{\bf E}}
\newcommand{\pr}{\mbox{\bf P}}
\newcommand{\E}{{\bf E}}
\newcommand{\R}{{\mathbb R}}
\newcommand{\one}{{\mathbf 1}}
\newcommand{\as}{\ \ \mbox{a.s.}}

\newcommand{\ep}{\varepsilon}

\newcommand\mnote[1]{} 
\newcommand\be{\begin{equation}}
\newcommand\bel[1]{{\mnote{#1}}\be\label{#1}}
\newcommand\ee{\end{equation}}
\newcommand\lb[1]{\label{#1}\mnote{#1}}

\title{Random walks that avoid their past convex hull}

\author{Omer Angel \and Itai Benjamini \and B\'alint Vir\'ag}

\begin{document}

\maketitle

\begin{abstract} We introduce planar random walk conditioned to avoid its
past convex hull, and we show that it escapes at a positive limsup speed.
Experimental results show that fluctuations from a limiting direction are
on the order of $n^{3/4}$.  This behavior is also observed for the
extremal investor, a natural financial model related to the planar walk.
\end{abstract}

\section{Introduction}

We consider the following random walk model and some closely
related models: inductively construct a sequence of points $x_i
\in \R^d$ by defining $x_0=0$ and $x_{n+1}$ to be uniformly
distributed on the sphere of radius 1 around $x_n$ but
conditioned so that the ``step'' segment
$\overline{x_nx_{n+1}}$ does not intersect the interior of the
convex hull of $\{x_0,\ldots,x_n\}$.


In the plane ($d=2$) this describes a frontier rancher who is
walking about and at each step increases his ranch by
``dragging'' with him the fence that defines it.


This model falls into the large category of self interacting
random walks, such as reinforced random walk or self-avoiding
walk. These models are difficult to analyze in general. The
reader should consult \cite{D}, \cite{L}, \cite{TW}, \cite{S},
and especially the survey papers \cite{T}, \cite{pemantle01}
for examples.

Over the next few sections we investigate the asymptotic
behavior of $x_n$. The main result, reported in Section
\ref{speed}, is that the planar rancher has positive lim sup
speed. We conjecture that the direction of $x_n$ converges a.s.
It would be natural to believe that the deviations of the
process from its eventual direction are diffusive (or they are
roughly described by a one-dimensional random walk with, say,
bounded increments). In Section \ref{simsec} we discuss
simulations indicating that this is not the case. Based on
these simulations we conjecture that at time $n$ the distance
of the farthest point on the path from the line $ox_n$ behaves
like $n^{3/4}$.

In Section \ref{investor} we study a related one-dimensional model that we
call the {\bf extremal investor}. This model describes what happens to the
value of a stock when the stockholder's decisions are influenced by best
and worst past performance in a simple way. Simulations for the critical
case of this process yield the same exponent $3/4$.

\section{{\bf Speed in 2 dimensions}}\label{speed}

Since the model provides some sort of ``repulsion'' of the
rancher from his past, it can be expected that the rancher will
escape faster than a regular random walk. In the 2 dimensional
case we show the following:

\begin{thm}\label{main}
If $d=2$ then there exists $s>0$ such that $\limsup \|x_n\|/n
\geq s$ a.s.
\end{thm}

This means that the rancher has positive lim sup speed. Our
simulations give an approximate speed of $0.314$.

The idea of the proof is to find a set of times of positive
density in which the expected gain in distance is bounded from
below.
\begin{figure}
\centering
\includegraphics[height=1.8in]{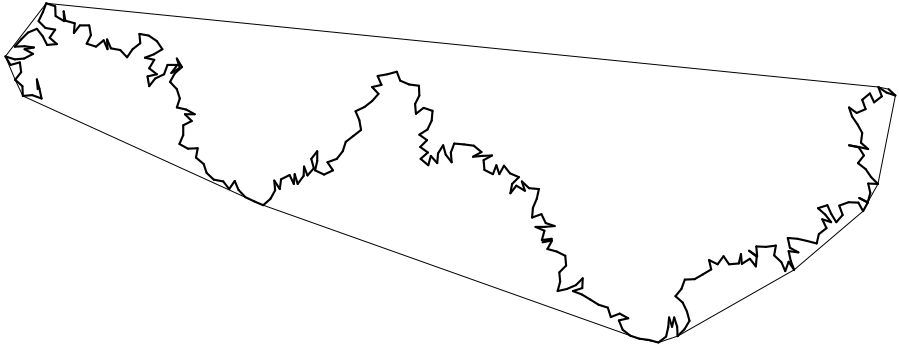}
\caption{300 steps of the rancher}
\end{figure}
There are two cases where the expected gain in distance can be
small. If, from the point of view of the rancher, the angle
that the ranch spans is very small then the next step is close
to uniform. The second problematic case is when the angle the
ranch spans is very close to $\pi$ with the direction of the
origin close to one of the ends. In this case the expected gain
in distance is also small.

\begin{proof}
Set $s_n = \|x_{n+1}\|-\|x_n\|$ and note that always $\E s_n
\geq 0$ since the legal directions of travel span an arc not
containing the origin.

Let $R_n$ denote the ranch at time $n$. If the angle of the
polygon $R_n$ at $x_n$ is in $[\ep,\pi-\ep]$,
then $\E s_n$ is bounded from below by some function of $\ep$.

If the angle is less then $\ep$ then we consider two
consecutive times. With probability at least half, the absolute
value of the angle $o x_n x_{n+1}$ is in $[\pi/4,3\pi/4]$. (Our
convention will be to regard an angle $xyz$ as a signed
quantity in $(-\pi,\pi]$.) In that case the angle the ranch
spans at time $n+1$ is large but not too large, and we have a
lower bound on $\E s_{n+1}$. If the first step is bad we just
use the bound $\E s_{n+1}>0$ and together we have a uniform
bound on $\E s_{n+1}$ in the case that the angle is small.

If the angle is large then we are in a tighter spot: it could
stay large for several steps. The rest of the proof consists of
showing that at positive fraction of time the angle is not
large.





We first introduce some notation. Consider the half-line
starting from $x_n$ that contains the edge of $R_n$ incident to
and clockwise from $x_n$. Let $y_n$ denote the intersection of
this half line and the boundary $C$ of the smallest disk about
the origin $o$ containing the ranch. Let $\alpha_n$ denote the
angle $\pi-ox_ny_n$, and let $\alpha'_n$ denote the analogous
angle in the counterclockwise direction. Let $d_n$ be the
distance between $C$ and $x_n$.

If $d_n$ is bounded above, then with probability bounded away
from 0, in a bounded number of steps the walk can get to a
position where $\ev s_{n'}>c>0.$ So it is suffices to show that
the Markov process $\{(R_n,x_n)\}$ returns to the set
$A=\{(R,x):\ d<d_*\}$ at a positive fraction of time.

To show this, we use a martingale argument; it suffices to
exhibit a non-negative function $f(R_n,x_n)$, so that the
expected increase in $f$ given the present is negative and
bounded away from zero when $(R_n,x_n)\not\in A$, and is
bounded from above when $(R_n,x_n)\in A$. The sufficiency of
the above is proved in Lemma \ref{technical} below; there take
$A_n$ to be the event $(R_n, x_n)\in A$, and $X_n=\|x_n\|$. We
now proceed to exhibit a function $f$ with the above properties.

The standard function that has this property is the expected
hitting time of $A$. We will try to guess this. The motivation
for our guess is the following heuristic picture. When the
angle $\alpha$ is small, it has a tendency to increase by a
quantity of order roughly $1/d$, and $d$ tends to decrease by a
quantity of order $\alpha$. This means that $d$ performs a
random walk with downward drift at least $1/d$, but this is not
enough for positive recurrence. So we have to wait for a few
steps for $\alpha$ to increase enough to provide sufficient
drift for $d$; the catch is that in every step $\alpha$ has a
chance of order $\alpha$ to decrease, and the same order of
chance to decrease to a fraction of its size. So $\alpha$ tends
to grow steadily and collapse suddenly. If the typical size is
$\alpha_*$, then it takes order $1/\alpha_*$ time to collapse.
During this time it grows by about $1/(d\alpha_*)$, which
should be on the order of the typical size $\alpha_*$, giving
$\alpha_*=d^{-1/2}$. This suggests that the process $d$ has
drift of this order, so the expected hitting time of $0$ is of
order $d^{3/2}$. A more accurate guess depends on $\alpha$,
too.

We define the functions $f_1(d)=d^{3/2}$, $f_2(d,\alpha)=-(
(cd^{1/2})\land(\alpha d))$, where $c=1/6$ is a constant, and
$f(d,\alpha,\alpha')=f_1(d)+ f_2(d,\alpha)+f_2(d,\alpha')$. It
is clear that $f(d_n,\alpha_n,\alpha'_n)$ can only increase by
a bounded amount on $A$. $f$ can be negative, but it is bounded
from below, which is sufficient. We want to show that given the
present outside $A$ the expected change in
$f(d_n,\alpha_n,\alpha'_n)$ is negative and bounded away from
zero.
\begin{figure}
\centering
\includegraphics[height=2.6in]{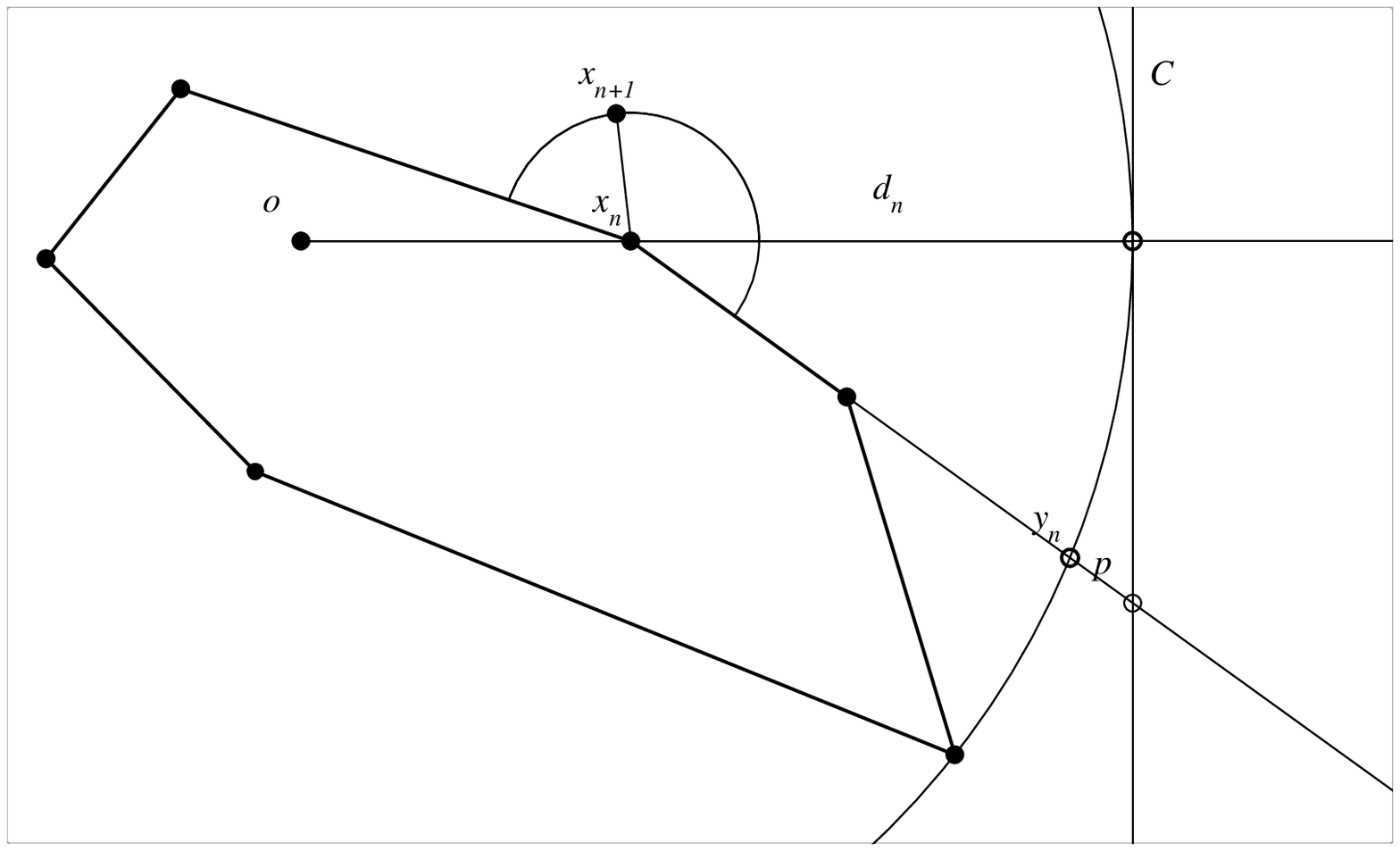}
\end{figure}
First we consider the expected change in $f_1$. All expected
values will be conditional on the information available at time
$n$.  To simplify notation, assume that the coordinates of
$x_n$ satisfy $x_{n,1}>0$, and $x_{n,2}=0$. We first bound the
expected decrease $d_n-d_{n+1}$.
 \bel{dnch}
 -\ev \Delta d 
  = \ev \|x_{n+1}\|-\|x_n\| \geq \ev x_{n+1,1}-x_{n,1}
\ee
 The right hand side can be computed directly. Let
$\beta=ox_nx_{n+1}-\pi$ denote the random angle of the $n$th
step. We keep our convention that $\beta \in (-\pi, \pi]$; for
example, $\beta=0$ means that the walker moved directly away
from $o$, $\beta>0$ means that the walker moved
``counterclockwise''. Then $\beta$ has uniform distribution on
$[-\alpha_n,\alpha'_n].$ We can then write the right of
(\ref{dnch}) as an integral
$$
 {1\over \alpha_n+\alpha_n'}
 \int_{-\alpha_n}^{\alpha_n'} \cos \beta \ d\beta \geq
 {\sin \alpha_n+\sin \alpha_n' \over 2\pi}
$$
 Using the fact that for $\Delta d$ bounded we have
 $$
 (d+\Delta d)^{3/2} = d^{3/2}+{3\over 2}d^{1/2}\Delta d + O(d^{-1/2})
 $$
 we bound
 \bel{df1}
 \ev \Delta f_1 \leq -{3 \over 4\pi}\left(\sin \alpha_n+ \sin \alpha'_n\right) d_n^{1/2}
  + \delta.
 \ee
Here and in the sequel $\delta$ denotes any quantity that
converges to $0$ if $d_*$ converges to $\infty$ ($d_*$ is a
constant to be set later so that $\delta$ is sufficiently
small).

We now proceed to bound the expected change in
$f_2(d_n,\alpha_n)$; denote this change by $\Delta f_2$. We
break up $\Delta f_2$ into important and unimportant parts:
 \begin{eqnarray*}
 \Delta f_2
 &=&  (cd_n^{1/2}\land \alpha_n d_n
    -cd_{n}^{1/2}\land \alpha_{n+1}d_n) \\
  &+&  (cd_n^{1/2}\land \alpha_{n+1}d_n
    -cd_{n+1}^{1/2}\land \alpha_{n+1}d_n)\\
  &+& (cd_{n+1}^{1/2}\land \alpha_{n+1}d_n
       -cd_{n+1}^{1/2}\land \alpha_{n+1}d_{n+1}).
 \end{eqnarray*}
The second term is bounded above by $c|d_{n+1}^{1/2}-
d_n^{1/2}|=\delta$, the third term is non-positive unless
$cd_{n+1}^{1/2}>\alpha_{n+1}d_n$, and then it can be at most
$\alpha_{n+1}|\Delta d| = \delta$. Thus important increase can
only come from the third term; call it $z$.  We examine three
cases according to the value of $\beta$.

{\it Event $B_2$:} $\beta\in[0,\pi-\alpha_n]$ (equivalently,
$x_{n+1}$ is on the side opposite of $R_n$ for the lines $ox_n$
and $x_ny_n$). Then
\begin{eqnarray}
\Delta\alpha &=& ox_ny_n - ox_{n+1}y_{n+1}
 \geq ox_ny_n - ox_{n+1}y_n \nonumber\\
 &=& x_nox_{n+1}+
x_{n+1}y_nx_n \geq x_{n+1} y_n x_n \geq 0. \label{dal1}
\end{eqnarray}
All inequalities follow from our assumption $B_2$. The equality
follows from the fact that the angles in the quadrangle
$ox_ny_nx_{n+1}$ add up to $2\pi$. A byproduct of (\ref{dal1})
is that $B_2$ implies $z\leq 0$. We now compute the last angle
in (\ref{dal1}) using a simple identity in the triangle
$x_ny_nx_{n+1}$:
 \bel{dal2}
 \dist(x_{n+1},y_n)  \sin (x_{n+1} y_n x_n) =
 \dist(x_{n},x_{n+1}) \sin (y_nx_nx_{n+1}) =
 \sin(\beta+\alpha_n).
 \ee
Now assume that $\alpha_n<(c-\delta)d_n^{-1/2}$. Then
$$
\dist(x_{n+1},y_n)\leq 1+\dist(x_n,y_n) \leq 1+\dist(x_n,p) =
1+d_n (\cos \alpha_n)^{-1} = d_n(1+\delta),
$$
where the point $p$ is the intersection of the tangent line to
$C$ at the ray $ox_n$ and the line $x_ny_n$.  We can then
conclude from (\ref{dal1}) and (\ref{dal2}) that
 \bel{dal3}
 \Delta \alpha \geq
  x_{n+1} y_n x_n
  \geq (1-\delta)\sin(\beta+\alpha_n)/d_n.
 \ee
The criterion $\alpha_n < (c-\delta)d_n^{-1/2}$ (for $\delta$
not too small) guarantees that the cutoff at $cd^{-1/2}$ does
not apply too early, and (\ref{dal3}) implies $z\leq
-(1-\delta)\sin(\beta+\alpha_n)$. Therefore
$$
 \ev [z;\;B_2] \leq -
 (\alpha_n+\alpha_n')^{-1}
 \int_0^{\pi-\alpha_n}
 \sin (\beta-\alpha_n) d\beta + \delta
 \leq -2\pi^{-1} + \delta.
$$

 {\it Event $B_3$:} $\beta>\pi-\alpha_n$. In this case
$R_{n+1}$ has an edge $x_{n+1} x_n$, and clearly
$\alpha_{n+1}>\pi-\beta$. Thus
$$\pr[0<z\mbox{ and }B_3]\leq\pr[\pi-\beta<cd_n^{-1/2}]\leq cd_n^{-1/2}\pi^{-1}.
$$

{\it Event $B_1$:} $\beta<0$. We can bound $\alpha_{n+1}$ below
by $\beta+\alpha_n$ as follows. First, note that
$\alpha_{n+1}=\pi-ox_{n+1}y_{n+1} \geq \pi-ox_{n+1}y_n$. Also
$\beta+\alpha_n=y_nx_nx_{n+1}=\pi-x_nx_{n+1}y_n-x_{n+1}y_nx_n$,
since the angles of a triangle add to $\pi$. We can split
$x_nx_{n+1}y_n=x_nx_{n+1}o + ox_{n+1}y_n$. Putting these
together we get $\alpha_{n+1}=\beta+\alpha_n+x_nx_{n+1}o +
x_{n+1}y_nx_n$, and since the latter two angles are small and
positive, $\alpha_{n+1}>\beta+\alpha_n$. Therefore
$$
 \pr[0<z\mbox{ and }B_1]\leq\pr[\beta+\alpha_n<cd_n^{-1/2}]\leq cd_n^{-1/2}\pi^{-1}.
$$
We now summarize our estimates. Since $z$ can be at most
$cd_n^{1/2}$, there is at most a bounded amount of positive
drift in $f_2$:
 $$
 \ev [\Delta f_2;\; B_1\cup B_3] \leq 2c^2\pi^{-1}
 +\delta.
 $$
If $\alpha_n<(c-\delta)d^{-1/2}$, then this is offset by the
negative drift
$$\ev[\Delta f_2;\;B_2] \leq - 2\pi^{-1} +\delta,
$$
and both inequalities have counterparts for $\alpha_n'$. If
$\alpha_n,\ \alpha_n'\geq (c-\delta)d^{-1/2}$, and at least one
of them is less than $\pi-\ep$, then from (\ref{df1}) we have
 $$
 \ev \Delta f_1\leq -3c/(4\pi)+\delta.
 $$
So for the cases covered so far,
$$
\ev \Delta f \leq (4c^2-2\land 3c/4)\pi^{-1}+\delta.
$$
This, for small $\delta$ and $c=1/6$, is negative and bounded
away from 0. The only remaining case is when
$\alpha_n,\alpha_n'>\pi-\ep$. We have seen that (looking at two
steps at a time)  $\ev d_{n+2}-d_n$ is bounded below by a fixed
constant, hence the expected decrease in $f_1$ is at least a
constant times $d_n^{1/2}$, which is enough to offset any
bounded positive drift in $f_2$. Another way to handle this
case is to add $f_3=\one(\alpha,\alpha'>\pi-\ep)$ to $f$.
\end{proof}

For the following lemma, we use the notation $\Delta_m a_n =
a_{n+m}-a_n$, and $\Delta a_n=\Delta_1 a_n$.

\begin{lemma}\label{technical}
Let $\{(X_n, f_n, A_n)\}$ be a sequence of triples adapted to
the increasing filtration $\{\mathcal F_n\}$ (with $\mathcal F_0$ trivial) 
so that  $X_n$,
$f_n$ are random variables and $A_n$ are events satisfying the
following. There exist positive constants $c_1, c_2, c_3, c_4$,
and a positive integer $m$, so we have a.s. for all $n$
\begin{eqnarray} 
|\Delta X_n| &\leq& 1, \nonumber \\
\ev [\Delta X_n\ |\ \mathcal F_n] &\geq& 0,\label{ketto} \\
\ev [\Delta_m X_n\ |\ \mathcal F_n, A_n] &>& c_1, \label{harom}\\
f_n &>& -c_2, \nonumber \\
\Delta f_n \one(A_n) &<& c_3, \label{ot}\\
\ev [\Delta f_n\ |\ \mathcal F_n, A_n^c] &<& -c_4. \label{hat}
\end{eqnarray}
Then for some positive constant $c_5$ we have
\begin{equation}\label{amikell}
  \limsup X_n/n
>c_5 \as
\end{equation}
\end{lemma}

\begin{proof} Let $G_n=\sum_{i=0}^{n-1} \one_{A_i}$, and let
$G_{n,k}=\sum_{i=0}^{n-1} \one_{A_{mi+k}}$, $0\leq k < n$. First we show 
that the $m+1$ processes
\begin{eqnarray}
 &\{c_1 G_{n,k} - X_{mn+k}\}_{n \geq 0}, \ \ \ \
0\leq k <m, \label{firstm}\\
&\{f_n - c_3 G_{n} + c_4 (n-G_{n})\}_{n \geq 0} \label{lastproc}
\end{eqnarray}
are supermartingales adapted to $\{\mathcal F_{mn+k}\}_{n\geq0}$, $0\leq k <m$, 
$\{\mathcal F_n\}_{n\geq 0}$,
respectively. For the first $m$  processes fix $k$, and note that 
$\ev[c_1(G_{n+1,k}-G_{n,k})\,|\,F_{mn+k}]=c_1\one(A_{mn+k})$. 
Consider
$$
\ev[c_1(G_{n+1,k}-G_{n,k})\,|\,\mathcal F_{mn+k}]+\ev[ -(X_{m(n+1)+k}-X_{mn+k})\,|\,\mathcal F_{mn+k}]
$$
If $A_{mn+k}$ happens, then the first term equals $c_1$, and the second
is less than $-c_1$ by (\ref{harom}). If $A_{mn+k}$ does not happen, 
then the first term equals $0$ and the second is nonpositive by
(\ref{ketto}).
Putting these two together shows that (\ref{firstm}) are supermartingales. 
For the last process, consider
$$
\ev[\Delta f_n\,|\, \mathcal F_n]+ \ev[-c_3 \Delta G_n\,|\,  \mathcal F_n]
+ \ev[c4(1-\Delta G_n)\,|\, \mathcal F_n].
$$
If $A_n$ happens, then the first term is less than $c_3$ by (\ref{ot}), 
the second term equals $-c_3$, and the last equals $0$. If $A_n$ does not 
happen, then 
the first term is less than $-c_4$ by (\ref{hat}), the second term equals $0$,
and the third equals $c_4$. In both cases we get that the process
(\ref{lastproc}) is a supermartingale. 

It follows from the supermartingale property that for some $c>0$ 
and all $n\geq 0$ we have
\begin{eqnarray}
 \ev X_{mn+k}&\geq& c_1 \ev G_{n,k} - c, \ \  \ \ \ 0\leq k <m,
 \label{xin}
 \\
 \ev G_n &\geq&
 c_4/(c_3+c_4) n -c.
 \label{gin}
\end{eqnarray}
 Since
$G_{mn}=G_{n,0}+ \ldots +G_{n,m-1}$, it follows from
(\ref{gin}) that for some $c_6>0$ and all large $n$ there is
$k=k(n)$, so that $\ev G_{n,k}> c_6 n$. Then for some $c_7>0$ we
have $\ev X_{nm+k}>c_7 n$ by (\ref{xin}). As a consequence, for
$Y_n=\max\{X_{mn}, \ldots, X_{mn+m-1}\}$ we have $\ev
Y_n> c_7 n$.

Thus for some $c_8<1$ we have $\ev(1-Y_n/(mn))<c_8$ for all
large $n$. Since $X_n\leq X_0+n$, we have
$Y_n\leq X_0+mn+m-1=mn+c_{9}$ and therefore $1-(Y_n-c_9)/(mn)\geq 0$.
Fatou's
lemma then implies $$
\ev \liminf (1-Y_n/(mn))=\ev \liminf (1-(Y_n-c_9)/(mn))\leq c_8,$$ 
 for some $c_{10}\in(c_8,1)$ 
Markov's inequality gives
$\pr(\liminf (1-Y_n/(mn))< c_8/c_{10})>1-c_{10}$. So
for some $c_5>0$,
$$\pr(\limsup X_n/n > c_5)>1-c_{10},$$
but we can repeat this argument while conditioning on the
$\sigma$-field $\mathcal F_t$ to get
$$
\pr(\limsup X_n/n > c_5 \ |\ \mathcal F_t)>1-c_{10}
$$
so letting $t\to \infty$ by L\'evy's 0-1 law we get
(\ref{amikell}).
\end{proof}

\section{Angular convergence, $d=2$}

\par
In the case $d=2$ we have seen that the rancher has positive
speed. This means he is similar to a random walk where the
radius is growing linearly and there is a random movement in
the angular direction. Since the distance is linear in $n$ we
have that the angular change is of order $n^{-1}$. If the signs
of the angular change were independent this would imply angular
convergence.

In our case the angular movements are positively correlated:
after a move in one direction the process tends to keep moving
in that direction. Simulations suggest that these correlations
are not enough to stop angular convergence, and we conjecture
that this is in fact the case.






%
\begin{figure}
\centering \label{sim}
\includegraphics[height=2.5in]{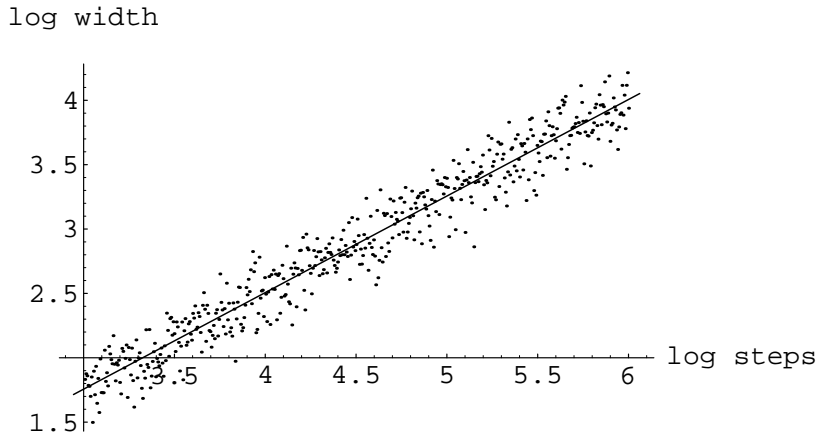}
\caption{The dimensions of the ranch}
\end{figure}
\section{Simulations and the exponent 3/4} \label{simsec}
 Computer simulations of the two dimensional process show
angular convergence to a random direction. We measured a
related quantity, the width $w_n$ of the path at time $n$,
defined as the distance of the farthest point on the path from
the line $ox_n$.

It is natural to guess that $w_n$ should behave as the maximum
of up to time $n$ of a one-dimensional Brownian motion, and
have a typical size of $n^{1/2}$. Our simulations, however,
show an entirely different picture. Figure \ref{sim} is a log
base $10$ plot of 500 realizations of $w_n$ on independent
processes. $n$ ranges from a thousand to a million steps
equally spaced on the log scale. The slope of the regression
line is $0.746$ (SE $0.008$). A regression line on the medians
of $1000$ measurements of walks of length $10^3,
10^4,10^5,10^6$ gave a value of $.75002$ (SE $0.002$). Based on
these simulations, we conjecture that $w_n$ behaves like
$n^{3/4}$. To put it rigorously in a weak form:

\begin{conj}\label{3/4}
For every $\ep>0$ we have $\pr[n^{3/4-\ep}<w_n<n^{3/4+\ep}]
\to 1$ as $n\to \infty$.
\end{conj}

\section{The extremal investor}\label{investor}

Stock or portfolio prices are often modeled by exponentiated
random walk or Brownian motion. In the simplest discrete-time
model, the log stock price, denoted $x_n$, changes every time
by an independent standard Gaussian random variable.

Ones decision whether to invest in, say, a mutual fund is often
based on past performance of the fund. Mutual fund companies
report past performance for periods ending at present; the
periods are often hand-picked to show the best possible
performance. The simplest such statistic is the overall best
performance over periods ending in the present. In terms of log
interest rate it is given by
 \bel{rate}
 r^{\max}_n=\max_{m<n}{x_n-x_m \over n-m },
 \ee
that is the maximal slope of lines intersecting the graph of
$x_n$ in both a past point and the present point.

A more cautious investor also looks at the worst performance
$r^{\min}_n$, given by (\ref{rate}) with a $\min$, and makes a
decision to buy, sell or hold accordingly, influencing the fund
price. In the simplest model, which we call the {\bf extremal
investor model}, the change in the log fund price given the
present is simply a Gaussian with standard deviation 1 and
expected value given by a fixed influence parameter $\alpha$
times the average of $r^{\max}$ and $r^{\min}$:
$$
x_{n+1} = x_n + \alpha {r^{\max}_n+r^{\min}_n\over 2} + \mbox{
standard Gaussian}.
$$
\begin{figure}
\centering \label{inv}
\includegraphics[height=2.2in]{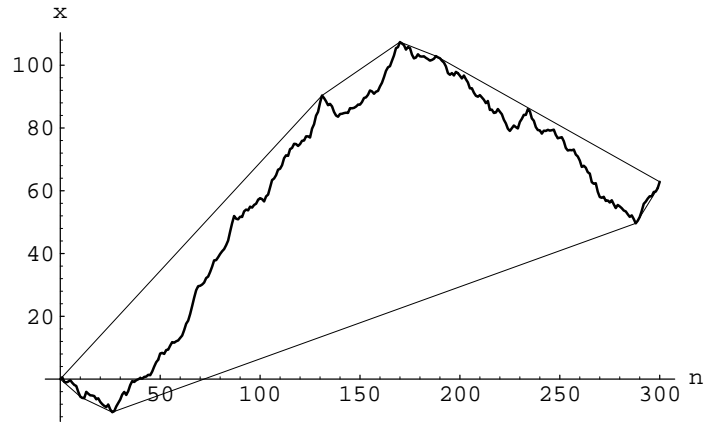}
\caption{The extremal investor process for $\alpha=1$ and its
enveloping curves}
\end{figure}

\noindent This process is related to the rancher in two dimensions, since
the future behavior of $x_n$ is influenced through the shape of
the convex hull of the graph of $x_n$ at the tip. Let $w_n$
denote the greatest distance between $x_n$ and the linear
interpolation from time zero to the present (assume $x_0=0$):
$$
w_n=\max_{m\leq n} \left| x_m-{m \over n}x_n \right|.
$$
We have the following version of Conjecture \ref{3/4}:
\begin{conj}\lb{3/4b}
 Let $\alpha=1$. For every $\ep>0$ we have
$\pr[n^{3/4-\ep}<w_n<n^{3/4+\ep}] \to 1$ as
$n\to \infty$.
\end{conj}
A moment of thought shows that for $\alpha>1$, $x_n$ will blow
up exponentially, so $\alpha_c=1$ is the critical parameter.
For $\alpha<1$ the behavior of $w_n$ seems to be governed by an
exponent between $1/2$ and $3/4$ depending on $\alpha$.
Simulations confirm Conjecture \ref{3/4b}. For $\alpha<1$ the
$x_n/n$ seems to converge to $0$, but in the case of
$\alpha=1$, it converges to a nontrivial random variable a.s.

\section{Questions}

\begin{ques}
Prove Theorem \ref{main} with $\liminf$ instead of $\limsup$.
\end{ques}

\begin{ques}
Show that for some $s>0$ the 2-dim rancher has speed $s$ a.s.
This could follow from some super-linearity result on the
rancher's travels.
\end{ques}

\begin{ques}
What is the behavior in higher dimensions? Does speed remain
positive? If not, is $\|x_n\| = O(\sqrt{n})$ or is it
significantly faster? What about convergence of direction?
\end{ques}

\begin{ques}
In $d=2$ what is the scaling limit of the path?
\end{ques}

\begin{ques}
If longer step sizes are allowed what happens when the tail is
thickened? Are there distributions which give positive speed
without convergence of direction?
\end{ques}

\end{document}